\def\ra{\rightarrow}
\def\ms{\mapsto}
\def\mod{\mathop{\;\rm mod}\nolimits}

\def\Hom{\mathop{\rm Hom}\nolimits}

\def\Z{\hbox{\bf Z}}
\def\Q{\hbox{\bf Q}}

\def\F{\hbox{\bf F}}

\newtheorem{theo}{Th\'eor\`eme}

\def\medpar{\par\medskip}

\def\thebibliography#1{\section*{Bibliographie\markboth
 {REFERENCES}{REFERENCES}}\list
 {[\arabic{enumi}]}{\settowidth\labelwidth{[#1]}\leftmargin\labelwidth
 \advance\leftmargin\labelsep
 \usecounter{enumi}}
 \def\newblock{\hskip .11em plus .33em minus -.07em}
 \sloppy
 \sfcode`\.=1000\relax}

\def\@begintheorem#1#2{\global\advance\@listdepth -1\relax
\list{}{\leftmargin 0pt\labelwidth -\labelsep}\item[{\sc #1\ #2}]\it}
\def\@endtheorem{\endlist\global\advance\@listdepth 1\relax}

\def\section{\@startsection {section}{1}{\z@}{3.5ex plus 1ex minus
.2ex}{2.3ex plus .2ex}{\large\bf}}
\def\subsection{\@startsection{subsection}{2}{\z@}{3.25ex plus 1ex minus
.2ex}{1.5ex plus .2ex}{\normalsize\bf}}
\def\subsubsection{\@startsection{subsubsection}{3}{\z@}{3.25ex plus 1ex
minus .2ex}{1.5ex plus .2ex}{\normalsize\bf}}

 \def\tableofcontents{\section*{Table des mati\`eres}
 \markboth{CONTENTS}{CONTENTS}
 \@starttoc{toc}}

\def\Cl{\rm Cl\;}
\documentstyle[A4,11pt]{article}
\begin{document}
G\'eom\'etrie alg\'ebrique/ {\it Algebraic geometry.}

\medpar
{\large Corps quadratiques dont le $5$-rang  du groupe des classes
est $\geq 3$.}

Jean-Fran\c cois Mestre.

{\bf R\'esum\'e.-} Nous prouvons qu'il existe une infinit\'e de corps
quadratiques r\'eels (resp. imaginaires) dont le  $5$-rang du groupe
des classes d'id\'eaux est $\geq 3$.

\medpar
{\large Quadratic fields whose $5$-rank is $\geq 3$.}

{\bf Abstract.-} We prove the existence of infinitely many real and
imaginary fields whose $5$-rank of the class group is $\geq 3$.
\vspace{5ex}

Soient $K$ un corps de nombres,  $\Cl K$ le groupe des classes d'id\'eaux de
$K$,
et $p$ un nombre premier. Par d\'efinition,
le $p$-rang du groupe des classes de $K$ est la dimension de $\Cl K/p\Cl K$
sur $\F_p$.
\medpar
Nous d\'emontrons ici le th\'eor\`eme suivant:

\begin{theo}
Il existe une infinit\'e de corps quadratiques r\'eels (resp.
imaginaires) dont le $5$-rang du groupe des classes est $\geq 3$.
\end{theo}

L'id\'ee de la d\'emonstration est la suivante: soit $E$ une courbe
elliptique d\'efinie sur $\Q$, munie d'un point $P$ d'ordre $5$ rationnel
sur $\Q$; si $F$ est la courbe quotient $E/<P>$, notons $\phi:\;\;E\ra F$
l'isog\'enie canonique de $E$ sur $F$.

Le lemme suivant, cons\'equence de \cite{RAYNAUD:schemas} et
du lemme $3$ de \cite{RAYNAUD:jacobienne}, m'a \'et\'e indiqu\'e par M.
Raynaud:

{\sc Lemme.-} {\it
 Supposons $E$ semi-stable
en tout nombre premier $p$; soient $K$ un corps quadratique, et
$O_K$ son anneau d'entiers. Si ${\cal E}$
(resp. ${\cal F}$) est le mod\`ele de N\'eron de $E$ (resp. $F$) sur $O_K$,
on a une suite exacte (de  sch\'emas en groupes sur $O_K$):
$$0\ra \Z/5\Z \ra {\cal E} \ra {\cal F}'\ra 0,$$
o\`u $\cal F'$ est un sous-sch\'ema en groupes ouvert de ${\cal F}$ contenant
la composante neutre ${\cal F}^0$ de ${\cal F}$.}
\medpar
Par suite, l'image r\'eciproque
par $\phi$ de tout point de ${\cal F}'(O_K)$ engendre une extension
ab\'elienne non ramifi\'ee de degr\'e divisant $5$ de $K$.

Soit  $$y^2+a_1xy+a_3y=x^3+a_2x^2+a_4x+a_6$$
une \'equation minimale de Weierstrass de $F$ sur $\Z$. Soit $Q=(x,y)$
un point de $F$ tel que $x\in \Q$; $Q$ est alors rationnel sur le corps
$K=\Q(y)$.
Notons $S$ l'ensemble fini des nombres premiers $p$
tels que le nombre de composantes connexes de la fibre en $p$
du mod\`ele de N\'eron de $F$  sur $\Z$ est divisible par $5$.
Supposons que, pour tout nombre premier
$p\in S$, $Q$ ne se r\'eduise pas $\mod p$ en le point singulier de
$F_{/\F_p}$. Alors le point $Q$ se prolonge en un
point de ${\cal F}'(O_K)$.
Par suite, si $L=K(\phi^{-1}(Q))$, $L$ est une extension non ramifi\'ee
de $K$. Le th\'eor\`eme d'irr\'eductibilit\'e de Hilbert permet
de montrer que, pour une infinit\'e de tels $x\in \Q$, $L/K$ est
de degr\'e $5$ (En fait, comme on le verra plus loin,
on peut donner des crit\`eres effectifs de congruence modulo des nombres
premiers convenables pour assurer que $L/K$ est de degr\'e $5$.)

Nous construisons dans la section suivante une courbe $X$, d\'efinie sur $\Q$,
poss\'edant les propri\'et\'es suivantes:

i) Il existe un rev\^etement $\psi$ de degr\'e $2$ de $X$ sur la droite
projective.

ii) Pour $1\leq i\leq 3$, il existe une courbe elliptique $E_i$, semi-stable
sur $\Z$, poss\'edant un point $\Q$-rationnel $P_i$, et un
rev\^etement ab\'elien $\tau_i$ de groupe de Galois $(\Z/2\Z)^2$ de $X$ sur
la courbe $F_i=E_i/<P_i>$.

Notons $\phi_i$ le morphisme canonique de $E_i$ sur $F_i$.
Nous montrons ensuite que, pour une infinit\'e de nombres rationnels $x$,
si $K=\Q(\psi^{-1}(x))$, les trois extensions
$L_i=K(\phi^{-1}(\tau_i(\psi^{-1}(x))))$ de $K$ sont ab\'eliennes de
degr\'e $5$, non ramifi\'ees et ind\'ependantes (i.e. les \'el\'ements
de $\Hom(\Cl K,\Z/5\Z)$ dont elles proviennent sont ind\'ependants), ce
qui prouve le th\'eor\`eme.

{\sc Remarque.- } Soit $E$ une courbe elliptique d\'efinie sur $\Q$,
munie d'un point rationnel $P$ d'ordre $n$, $n$  entier $\geq 1$. Soit $F$
la courbe quotient $E/<P>$, $\phi:\;\;E\ra F$ le
morphisme canonique, et $y^2=f(x)$ un mod\`ele de Weierstrass
de $F$. Soit $x\in \Q$, et $Q$ l'un des deux points de $F$ d'abscisse $x$.
Si $K=\Q(\sqrt{f(x)})$,
le th\'eor\`eme de Chevalley-Weil permet de montrer que,
d\`es que la valuation de $x$ est suffisamment n\'egative en chaque
nombre premier o\`u $E$ a mauvaise r\'eduction, l'extension $K(\phi^{-1}(Q))/K$
est non ramifi\'ee. L'hypoth\`ese ``{\it $E$ semi-stable}'' n'est donc
pas n\'ecessaire. N\'eanmoins, dans le cas o\`u elle n'est pas v\'erifi\'ee,
les conditions de congruence sur $x$ sont plus d\'elicates \`a d\'eterminer.

\section{Construction de la courbe $X$}
\subsection{Construction de $E$ et $F$}
La courbe modulaire $X_1(10)$, classifiant les courbes elliptiques
munies d'un point d'ordre $10$, est de genre $0$, et a \'et\'e param\'etr\'ee
par Kubert \cite{KUBERT:universal}: il construit une courbe elliptique $E$
rationnelle sur $\Q(f)$, o\`u $f$ est un param\`etre,
poss\'edant un point $P_0$ d'ordre $10$
rationnel sur $\Q(f)$.

Si l'on pose $f=(u+1)/2$, et apr\`es un changement de variables, on
trouve comme \'equation de $E$:
$$y^2=(x^2-u(u^2+u-1))(8xu^2+(u^2+1)(u^4-2u^3-6u^2+2u+1)) .$$
Le point $P_0$ d'ordre $10$ est le point d'abscisse
$-{\frac {\left (u^{4}-2\,u^{3}-6\,u^{2}+2\,u+1\right )\left (u^{2}+1
\right )}{8\,u^{2}}}.$

Les formules de V\'elu \cite{VELU:isog} permettent
alors d'obtenir une \'equation de la courbe $F$ quotient de $E$ par
le groupe d'ordre $5$ engendr\'e par $2P_0$;
une \'equation de $F$ est donn\'ee par
$y^2=g_u(x)$, o\`u $$g_u(x)=(x^2-u(u^2+u-1))h_u(x)\;\;\;
{\rm et}\;\;\; h_u(x)=8(u^2+u-1)^2x+(u^2+1)(u^4+22u^3-6u^2-22u+1).$$

De plus, si $u\in \Q$, la condition $u\equiv \pm 1 \mod 5$ assure que les
courbes $E$ et $F$ sont semi-stables sur $\Z$.
\subsection{Construction de $X$}
Si
$$\left\{\begin{array}{l}
u_1=(t^2+t-1)/(t^2+t+1),\\ u_2=-(t^2+3t+1)/(t^2+t+1),\\
u_3=-(t^2-t-1)/(t^2+t+1)
\end{array}\right.$$
on a
$$u_1(u_1^2+u_1-1)=u_2(u_2^2+u_2-1)=u_3(u_3^2+u_3-1).$$

La courbe $X$, d\'efinie sur $\Q(t)$,  normalis\'ee de la courbe d'\'equations
$$\left\{\begin{array}{l}
y_1^2=g_{u_1}(x),\\ y_2^2=g_{u_2}(x),\\y_3^2=g_{u_3}(x)
\end{array}\right.$$
est
de genre $5$.
L'application
$\phi:\;\;(x,y_1,y_2,y_3)\ms (x,v=y_1/y_2,w=y_1/y_3)$
de $X$ sur la courbe $C$ de genre $0$ et
d'\'equations $$h_{u_1}(x)=v^2h_{u_2}(x)=w^2h_{u_3}(x)$$
est de degr\'e $2$; les quatre points de $C$ de coordonn\'ees $(x,v,w)=(\infty,
\pm u_2/u_1,\pm u_3/u_1)$ sont rationnels sur $\Q(t)$
et sont des points de ramification de $\phi$,
donc $C$ est $\Q(t)$-isomorphe \`a la
droite projective, et $X$ est hyperelliptique, rev\^etement double de la
droite projective, et poss\`ede quatre points de Weierstrass rationnels
sur $\Q(t)$.

Sp\'ecialisons en $t=4$ les formules de la section pr\'ec\'edente.
On trouve $u_1=19/21$, $u_2=-29/21$ et $u_3=-11/21$. Les courbes
$E_i$ et $F_i$, pour $1\leq i\leq 3$, sont donc semi-stables sur $\Z$.

Une param\'etrisation de la courbe $C$ est donn\'ee
par

$$v=\frac{29}{19}\frac{53719189282 z^2+26766692861}{ 53719189282 z^2
- 283246634396 z- 26766692861},$$
$$w=\frac{11}{19}\frac{53719189282 z^2+26766692861}{
53719189282 z^2+ 20305766998 z- 26766692861},$$

$$x=\frac{c_4z^4+c_3z^3+c_2z^2+c_1z+c_0}{5167944494559 (4883562662 z+922989409)
(11 z -29) z},$$
avec
$$\begin{array}{ll}
c_0=343898806423252015354080,&c_1=- 411804539876837130626339,\\
c_2=- 642297925780193483509181,&c_3=826467660375890872281118,\\
c_4=1385160622615364964251520.&\end{array}$$

On obtient alors une \'equation hyperelliptique de $X$ en substituant
la fraction rationnelle $x(z)$ ci-dessus dans, par exemple, l'\'equation
de $E_1$, \`a savoir
$$y^2=42(44876601 x-133597561)(9261 x^2-6061).$$
\`A toute valeur rationnelle de $z$
distincte des p\^oles de $x(z)$, on associe ainsi le corps $K=\Q(y)$.
\subsection{Conditions sur $z$ pour que le $5$-rang de $\Cl K$ soit
$\geq 3$}
Si $p$ est un nombre premier, notons $v_p$ la valuation $p$-adique.
Le calcul montre qu'un point de $F_1(K)$
(resp. $F_2(K)$, resp. $F_3(K)$) se prolonge \`a ${\cal F_1}'(O_K)$
(resp. ${\cal F_2}'(O_K)$, resp. ${\cal F_3}'(O_K)$) si et seulement si
son abscisse $x$ v\'erifie $v_{11}(x)\leq -2$, $v_{29}(x)\leq -2$ et
$x\not\equiv 77 \mod 419$ (resp. $v_{11}(x)\leq -2$, $v_{19}(x)\leq -2$,
$x\not\equiv 677 \mod 709$, resp. $v_{19}(x)\leq -2$, $v_{29}(x)\leq -2$,
$x\not\equiv 36 \mod 151$).
Si
\begin{equation}
z\equiv 0 \mod 11.19.29\;\;\;{\rm et}\;\;\;
z\not\equiv \pm 86\mod 419,\end{equation}
les conditions de congruence ci-dessus
sont remplies, et les points correspondants de $F_i(K)$, $i=1,2,3$, se
prolongent en des points de $F'_i(O_K)$.

De plus, soit  $l_1=163$, $l_2=701$ et $l_3=1277$;  supposons
\begin{equation}
z\equiv 1\mod l_1l_2l_3.\end{equation} Alors:

i) Les id\'eaux $(l_i)$, $i=1,2,3$, se
d\'ecomposent chacun dans $K$ en deux id\'eaux ${\cal P_i}\overline{\cal
P_i}$,

ii) ${\cal P_1}$ (resp. ${\cal P_2}$, resp. ${\cal P_3}$)
est d\'ecompos\'e (resp.
d\'ecompos\'e, resp. inerte) dans $L_1$,

iii) ${\cal P_1}$ (resp. ${\cal P_2}$, resp. ${\cal P_3}$)
est inerte (resp. d\'ecompos\'e, resp.
 d\'ecompos\'e) dans $L_2$,

 iv) ${\cal P_1}$ (resp. ${\cal P_2}$, resp. ${\cal P_3}$)
 est d\'ecompos\'e (resp. inerte,
 resp. d\'ecompos\'e) dans $L_3$.

 Ceci assure que les extensions $L_i$, $i=1,2,3$, sont ind\'ependantes.

 Par suite, d\`es que $z$ v\'erifie les congruences $(1)$ et $(2)$ ci-dessus,
 le $5$-rang du groupe des classes de $K$ est $\geq 3$.

Le fait que, lorsque $z$ parcourt une infinit\'e de valeurs
rationnelles, on obtient une infinit\'e de corps quadratiques $K$
provient par exemple du th\'eor\`eme de Faltings (la courbe $X$,
\'etant de genre $>1$, n'a qu'un nombre fini de points rationnels
dans un corps $K$ fix\'e.)

 De plus, comme $X$ a trois points de Weierstrass rationnels, donc
 r\'eels, on obtient ainsi une infinit\'e de corps quadratiques r\'eels
 (resp. imaginaires) dont le $5$-rang du groupe des classes est
 $\geq 3$. Par exemple, si $z$ est $>0$ (resp. $<0$) et suffisamment proche
 de $0$
 (pour la topologie usuelle), $K$ est quadratique imaginaire (resp.
 r\'eel).

\end{document}